\documentclass[12pt,thmsa]{article}
\usepackage{amsfonts}
\usepackage{amsmath}

\input tcilatex
\begin{document}

\author{Nick Laskin\thanks{\textit{E-mail address}: nlaskin@rocketmail.com}}
\title{\textbf{A New Family of Fractional Counting Probability Distributions}%
}
\date{TopQuark Inc.\\
Toronto, ON, M6P 2P2\\
Canada}
\maketitle

\begin{abstract}
A new family of fractional counting processes based on a three-parameter
generalized Mittag-Leffler function was introduced and studied.

As applications we develop a fractional generalized compound process,
introduce and develop fractional generalized Bell polynomials and numbers,
fractional generalized Stirling numbers of the second kind, and a new family
of quantum coherent states. Stretched quantum coherent states, which are a
generalization of the famous Schr\"{o}dinger-Glauber coherent states, were
also introduced and studied.

In particular cases, the presented results reproduce known equations for
Poisson and fractional Poisson probability distributions, Bell numbers and
fractional Bell numbers, Stirling numbers and fractional Stirling numbers of
the second kind, as well as for known quantum coherent states.


\end{abstract}

\section{Introduction}

To study counting processes with a long memory, the fractional Poisson
process was invented in \cite{Laskin1}. The fractional Poisson probability
distribution captures long memory phenomena which manifests itself in the
non-exponential waiting time probability distribution function, observed
empirically in complex quantum and classical systems. An example of a
quantum system is the fluorescence intermittency of CdSe quantum dots, that
is, the fluorescence of single nanocrystals exhibits intermittent behavior,
namely a sequence of \textquotedblleft light on\textquotedblright\ and
\textquotedblleft light off\textquotedblright\ states deviate from Poisson
statistics \cite{Kuno}. As examples of classical systems, we mention the
distribution of waiting times between two consecutive transactions in
financial markets \cite{Sabatelli} and another one that comes from network
communication systems, where the durations of network sessions or
connections exhibit non-exponential behavior \cite{Willinger}.

Here, we introduce and develop a new family of fractional generalized
counting processes and probability distribution functions based on a
three-parameter generalized Mittag-Leffler function, known as the Prabhakar
function \cite{Prabhakar}. The fundamental prerequisite for this is the fact
that the three-parameter generalized Mittag-Leffler function is a completely
monotonic function. A new family of fractional generalized counting
probability distribution functions includes, as special cases, the famous
Poisson probability distribution and the fractional Poisson probability
distribution. We present the foundations and applications of the new
fractional generalized counting probability distribution.

As applications, we develop a fractional generalized compound process,
fractional generalized Bell polynomials and numbers, fractional generalized
Stirling numbers of the second kind and new quantum coherent states. As a
particular member of the new coherent states, the stretched quantum coherent
states, which are a generalization of the famous Schr\"{o}dinger-Glauber
coherent states \cite{Schrodinger} - \cite{Klauder}, were introduced and
studied.

With certain sets of fractality parameters, all our new results fit into the
well-known equations for the Poisson process and probability distribution,
the fractional Poisson process and probability distribution, and famous Schr%
\"{o}dinger-Glauber coherent states.

The paper is organized as follows.

In Section 1, a new family of fractional generalized counting processes was
introduced. In Section 2, the probability and moment generating functions of
the newly introduced probability distribution were calculated. The equations
for the mean and variance were obtained in Section 3. Waiting time
distribution for the fractional generalized counting process was found in
Section 4. Applications of the family of fractional generalized processes
and probability distributions are presented in Section 5. Applications
include the fractional generalized compound process, fractional generalized
Bell polynomials and numbers, fractional generalized Stirling numbers of the
second kind and new quantum coherent states. Conclusion summarizes our
findings.

\section{A new family of fractional counting probability distribution
functions}

Based on the developments presented in \cite{Laskin1} and \cite{Laskin2}, we
introduce the following new \textit{fractional generalized counting process}
whose probability distribution function is defined as

\begin{equation}
P_{\mu ,\nu }^{\gamma ,\sigma }(n,t)=\Gamma (\nu )\frac{(-z)^{n}}{n!}\frac{%
d^{n}}{dz^{n}}\left( E_{\mu ,\nu }^{\gamma }(z)\right) |_{z=-\lambda
_{\sigma }t^{\sigma }},\qquad t\geq 0,  \label{eq1.1}
\end{equation}%
\begin{equation*}
P_{\mu ,\nu }^{\gamma ,\sigma }(0,t)=\Gamma (\nu )E_{\mu ,\nu }^{\gamma
}(-\lambda _{\sigma }t^{\sigma }),
\end{equation*}

where $E_{\mu ,\nu }^{\gamma }(z)$ is the generalized three-parameter
Mittag-Leffler function given by its power series \cite{Prabhakar}

\begin{equation}
E_{\mu ,\nu }^{\gamma }(z)=\sum\limits_{m=0}^{\infty }\frac{(\gamma )_{m}}{%
m!\Gamma (\mu m+\nu )}z^{m},  \label{eq1.2}
\end{equation}

with $(\gamma )_{m}$ defined as

\begin{equation}
(\gamma )_{m}=\frac{\Gamma (\gamma +m)}{\Gamma (\gamma )},  \label{eq1.3}
\end{equation}

here $\Gamma (\alpha )$ is the Gamma function, which has the familiar
representation

\begin{equation*}
\Gamma (\alpha )=\int\limits_{0}^{\infty }dte^{-t}t^{\alpha -1},
\end{equation*}

the fractality parameters $\mu $, $\nu ,$\ $\gamma $ and $\sigma $ satisfy
the conditions%
\begin{equation}
0<\mu \leq 1,\quad \gamma >0,\quad \nu \geq \mu \gamma \quad \text{and\quad }%
0<\sigma \leq 1,  \label{eq1.4}
\end{equation}

and parameter $\lambda _{\sigma }>0$. Parameters $\mu $, $\nu ,$ $\gamma $
and $\sigma $ are dimensionless, and $\lambda _{\sigma }$ has units of [time]%
$^{-\sigma }$.

Note that the function $P_{\mu ,\nu }^{\gamma ,\sigma }(n,t)$, introduced by
Eq.(\ref{eq1.1}), is a probability distribution function. Indeed, $P_{\mu
,\nu }^{\gamma ,\sigma }(n,t)$ is normalized

\begin{equation}
\sum\limits_{n=0}^{\infty }P_{\mu ,\nu }^{\gamma ,\sigma }(n,t)=\Gamma (\nu
)E_{\mu ,\nu }^{\gamma }(0)=1,  \label{eq1.5}
\end{equation}

and positive due to the complete monotonicity of the three-parameter
generalized Mittag-Leffler function $E_{\mu ,\nu }^{\gamma }(-x)$ \cite%
{Gorska},

\begin{equation}
(-1)^{n}\frac{d^{n}E_{\mu ,\nu }^{\gamma }(-x)}{dx^{n}}\geq 0,\qquad x\geq 0,
\label{eq1.6}
\end{equation}

with parameters $\mu $, $\nu $\ and $\gamma $ satisfying the conditions
given by Eq.(\ref{eq1.4}).

Given Eq.(\ref{eq1.1}) and replacing $\lambda _{\sigma }t^{\sigma }$ with $%
x^{\sigma }$, we introduce a new family of fractional generalized
probability distributions $P_{\mu ,\nu }^{\gamma ,\sigma }(x)$,

\begin{equation}
P_{\mu ,\nu }^{\gamma ,\sigma }(n,x)=\Gamma (\nu )\frac{(-z)^{n}}{n!}\frac{%
d^{n}}{dz^{n}}E_{\mu ,\nu }^{\gamma }(z)|_{z=-x^{\sigma }},\qquad x\geq 0,
\label{eq1.6a}
\end{equation}

where all notations are defined above.

Equation (\ref{eq1.1}) is a generalization of formula (30) from \cite%
{Laskin1}, which is reproduced\footnote{%
Note that the parameter $\nu $ in formula (30) from \cite{Laskin1} has units
of [time]$^{-\mu }$. In our notations, this parameter is equal to $\lambda
_{\sigma }$ in Eq.(\ref{eq1.1}) when $\sigma =\mu $.} in the case when $%
0<\mu \leq 1$, $\nu =\gamma =1$ and $\sigma =\mu $.

The function $P_{\mu ,\nu }^{\gamma ,\sigma }(n,t)$ introduced by Eq.(\ref%
{eq1.1}) and function $P_{\mu ,\nu }^{\gamma ,\sigma }(n,x)$ defined by Eq.(%
\ref{eq1.6a}) describe a family of new probability distributions of counting
processes and a family of new fractional generalized probability
distributions, respectively. The family includes the probability
distribution function $P(n,t)$ of the well-known Poisson process, when $\mu
=\nu =\gamma =\sigma =1$,

\begin{equation}
P(n,t)\equiv P_{1,1}^{1,1}(n,t)=\Gamma (1)\frac{(-z)^{n}}{n!}\frac{d^{n}}{%
dz^{n}}\left( E_{1,1}^{1}(z)\right) |_{z=-\lambda _{1}t}=  \label{eq1.7}
\end{equation}%
\begin{equation*}
=\frac{(-z)^{n}}{n!}\frac{d^{n}}{dz^{n}}(e^{z})|_{z=-\lambda _{1}t}=\frac{%
(\lambda _{1}t)^{n}}{n!}e^{-\lambda _{1}t},
\end{equation*}

and the Poisson probability distribution $P(n,x)$,

\begin{equation}
P(n,x)\equiv P_{1,1}^{1,1}(n,x)=\frac{x^{n}}{n!}e^{-x},  \label{eq1.7a}
\end{equation}

where we took into account that

\begin{equation*}
E_{1,1}^{1}(z)=e^{z},
\end{equation*}

as it follows from Eq.(\ref{eq1.2}).

When $0<\mu \leq 1$, $\nu =\gamma =1$ and $\sigma =\mu $ the family gives us
the fractional Poisson process $P_{\mu }(n,t)$ \cite{Laskin1}

\begin{equation}
P_{\mu }(n,t)\equiv P_{\mu ,1}^{1,\mu }(n,t)=\frac{(-z)^{n}}{n!}\frac{d^{n}}{%
dz^{n}}\left( E_{\mu ,1}^{1}(z)\right) |_{z=-\lambda _{\mu }t^{\mu }},
\label{eq1.8}
\end{equation}

or in the power series form

\begin{equation}
P_{\mu }(n,t)=\frac{(\lambda _{\mu }t^{\mu })^{n}}{n!}\sum\limits_{k=0}^{%
\infty }\frac{(k+n)!}{k!}\frac{(-\lambda _{\mu }t^{\mu })^{k}}{\Gamma (\mu
(k+n)+1)},  \label{eq1.10}
\end{equation}

and the fractional Poisson probability distribution $P_{\mu }(n,x)$,

\begin{equation}
P_{\mu }(n,x)=\frac{x^{\mu n}}{n!}\sum\limits_{k=0}^{\infty }\frac{(k+n)!}{k!%
}\frac{(-x^{\mu })^{k}}{\Gamma (\mu (k+n)+1)}.  \label{eq1.10_1}
\end{equation}

It is easy to see that for $\mu =1$ Eqs.(\ref{eq1.10}) and (\ref{eq1.10_1})
are transformed into Eqs.(\ref{eq1.7}) and (\ref{eq1.7a}).

For $\mu =$ $\nu =\gamma =1$ and $0<\sigma \leq 1$ we obtain the process
whose probability distribution function $P_{\sigma }(n,t)$ is

\begin{equation}
P_{\sigma }(n,t)=P_{1,1}^{1,\sigma }(n,t)=\frac{(\lambda _{\sigma }t^{\sigma
})^{n}}{n!}\exp (-\lambda _{\sigma }t^{\sigma }),  \label{eq1.10_a}
\end{equation}

and the probability distribution $P_{\sigma }(n,x)$,

\begin{equation}
P_{\sigma }(n,x)=P_{1,1}^{1,\sigma }(n,x)=\frac{x^{\sigma n}}{n!}\exp
(-x^{\sigma }),\qquad x\geq 0.  \label{eq1.10b}
\end{equation}

Note that $P_{\sigma }(n,x)$ is normalized and positive due to the complete
monotonicity of function $\exp (-x^{\sigma })$ for $0<\sigma \leq 1$.

The process Eqs.(\ref{eq1.10_a}) can be considered as a non-homogeneous
Poisson process whose \textquotedblleft rate\textquotedblright\ $r(t)$
changes with time,

\begin{equation}
r(t)=\sigma \lambda _{\sigma }t^{\sigma -1}.  \label{eq.1.10c}
\end{equation}

Then $P_{\sigma }(n,t)$ can be represented as follows,

\begin{equation}
P_{\sigma }(n,t)=\frac{(\Lambda _{\sigma }(t))^{n}}{n!}\exp (-\Lambda
_{\sigma }(t)),  \label{eq1.10d}
\end{equation}

where

\begin{equation}
\Lambda _{\sigma }(t)=\int\limits_{0}^{t}d\tau r(\tau )=\sigma \lambda
_{\sigma }\int\limits_{0}^{t}d\tau \tau ^{\sigma -1}=\lambda _{\sigma
}t^{\sigma }.  \label{eq.1.10e}
\end{equation}

Finally, let's find an alternative representation of the probability
distribution $P_{\mu ,\nu }^{\gamma ,\sigma }(n,t)$ of the fractional
generalized counting process introduced by Eq.(\ref{eq1.1}). 

Calculation of the first order derivative from $E_{\mu ,\nu }^{\gamma }(z)$
defined by Eq.(\ref{eq1.2}) yields

\begin{equation}
\frac{d}{dz}\left( E_{\mu ,\nu }^{\gamma }(z)\right) =\gamma E_{\mu ,\nu
+\mu }^{\gamma +1}(z).  \label{eq1.11}
\end{equation}

Then iterating Eq.(\ref{eq1.11}) gives us \cite{Prabhakar}

\begin{equation}
\frac{d^{m}}{dz^{m}}\left( E_{\mu ,\nu }^{\gamma }(z)\right) =\gamma (\gamma
+1)...(\gamma +m-1)E_{\mu ,\nu +m\mu }^{\gamma +m}(z)=(\gamma )_{m}E_{\mu
,\nu +m\mu }^{\gamma +m}(z),  \label{eq1.12}
\end{equation}

for any $m\in 
\mathbb{N}
$. Thus, with the help of Eq.(\ref{eq1.12}) the probability distribution $%
P_{\mu ,\nu }^{\gamma ,\sigma }(n,t)$ of the fractional generalized counting
process can be expressed as follows,

\begin{equation}
P_{\mu ,\nu }^{\gamma ,\sigma }(n,t)=\Gamma (\nu )\frac{(-z)^{n}}{n!}(\gamma
)_{n}E_{\mu ,\nu +n\mu }^{\gamma +n}(z)|_{z=-\lambda _{\mu }t^{\mu }}=
\label{eq1.13}
\end{equation}

\begin{equation*}
\Gamma (\nu )\frac{(\lambda _{\sigma }t^{\sigma })^{n}}{n!}(\gamma
)_{n}E_{\mu ,\nu +n\mu }^{\gamma +n}(-\lambda _{\sigma }t^{\sigma }),
\end{equation*}

or, using a power series

\begin{equation}
P_{\mu ,\nu }^{\gamma ,\sigma }(n,t)=(\gamma )_{n}\Gamma (\nu )\frac{%
(\lambda _{\sigma }t^{\sigma })^{n}}{n!}\sum\limits_{k=0}^{\infty }\frac{%
(\gamma +n)_{k}}{k!\Gamma (\mu (n+k)+\nu )}(-\lambda _{\sigma }t^{\sigma
})^{k},  \label{eq1.14}
\end{equation}

where the parameters $\mu $, $\nu ,$\ $\gamma $ and $\sigma $ satisfy the
condition given by Eq.(\ref{eq1.4}).

For a fractional generalized probability distribution $P_{\mu ,\nu }^{\gamma
,\sigma }(n,x)$ we have the following power series

\begin{equation}
P_{\mu ,\nu }^{\gamma ,\sigma }(n,x)=(\gamma )_{n}\Gamma (\nu )\frac{%
(x^{\sigma })^{n}}{n!}\sum\limits_{k=0}^{\infty }\frac{(\gamma +n)_{k}}{%
k!\Gamma (\mu (n+k)+\nu )}(-x^{\sigma })^{k},\qquad x\geq 0.  \label{eq1.15}
\end{equation}

The function $P_{\mu ,\nu }^{\gamma ,\sigma }(n,t)$ defined by Eqs.(\ref%
{eq1.1}) and (\ref{eq1.14}) gives us the probability that in the time
interval $[0,t]$ we observe $n$ counts driven by a time-fractional stream.
The definitions given by Eqs.(\ref{eq1.1}) and (\ref{eq1.14}) introduce a
new family of fractional generalized counting processes that includes two
well-known counting processes:\ the Poisson process and the fractional
Poisson process.

The definitions given by Eqs.(\ref{eq1.6a}) and (\ref{eq1.15}) introduce a
new family of fractional generalized probability distributions that includes
two well-known probability distributions: the Poisson distribution and the
fractional Poisson distribution.

Some other particular values of the parameters $\mu ,\nu ,\gamma $ and $%
\sigma $ bring us new members of the family of fractional generalized
counting processes and fractional generalized probability distributions. For
example, when $0<\mu \leq 1$, $\gamma =1$, $\nu \geq \mu $ and $0<\sigma
\leq 1$ we obtain a fractional generalized counting process whose
probability distribution $P_{\mu ,\nu }^{1,\sigma }(n,t)$ is equal to

\begin{equation}
P_{\mu ,\nu }^{1,\sigma }(n,t)=\Gamma (\nu )\frac{(\lambda _{\sigma
}t^{\sigma })^{n}}{n!}\sum\limits_{k=0}^{\infty }\frac{(n+k)!}{k!\Gamma (\mu
(n+k)+\nu )}(-\lambda _{\sigma }t^{\sigma })^{k}  \label{eq1.15a}
\end{equation}

and probability distribution function $P_{\mu ,\nu }^{1,\sigma }(n,x)$ is

\begin{equation}
P_{\mu ,\nu }^{1,\sigma }(n,x)=\Gamma (\nu )\frac{(x^{\sigma })^{n}}{n!}%
\sum\limits_{k=0}^{\infty }\frac{(n+k)!}{k!\Gamma (\mu (n+k)+\nu )}%
(-x^{\sigma })^{k}.  \label{eq1.15b}
\end{equation}

When $\mu =1$, $\gamma >0$, $\nu \geq \gamma $ and $0<\sigma \leq 1$ we
obtain the fractional generalized counting process whose probability
distribution $P_{1,\nu }^{\gamma ,\sigma }(n,t)$ is equal to

\begin{equation}
P_{1,\nu }^{\gamma ,\sigma }(n,t)=(\gamma )_{n}\Gamma (\nu )\frac{(\lambda
_{\sigma }t^{\sigma })^{n}}{n!}\sum\limits_{k=0}^{\infty }\frac{(\gamma
+n)_{k}}{k!\Gamma (n+k+\nu )}(-\lambda _{\sigma }t^{\sigma })^{k},
\label{eq1.16}
\end{equation}

and probability distribution function $P_{1,\nu }^{\gamma ,\sigma }(n,x)$ is

\begin{equation}
P_{1,\nu }^{\gamma ,\sigma }(n,x)=(\gamma )_{n}\Gamma (\nu )\frac{(x^{\sigma
})^{n}}{n!}\sum\limits_{k=0}^{\infty }\frac{(\gamma +n)_{k}}{k!\Gamma
(n+k+\nu )}(-x^{\sigma })^{k}.  \label{eq1.16a}
\end{equation}

\section{Generating functions}

\subsection{Probability generating function}

The probability generating function $G_{\mu ,\nu }^{\gamma ,\sigma }(s,t)$
is defined as

\begin{equation}
G_{\mu ,\nu }^{\gamma ,\sigma }(s,t)=\sum\limits_{n=0}^{\infty }s^{n}P_{\mu
,\nu }^{\gamma ,\sigma }(n,t),  \label{eq2.1}
\end{equation}

where $P_{\mu ,\nu }^{\gamma ,\sigma }(n,t)$ is given by Eq.(\ref{eq1.1}).
In terms of $G_{\mu ,\nu }^{\gamma ,\sigma }(s,t)$ the probability
distribution function is

\begin{equation}
P_{\mu ,\nu }^{\gamma ,\sigma }(n,t)=\frac{1}{n!}\frac{\partial ^{n}}{%
\partial s^{n}}G_{\mu ,\nu }^{\gamma ,\sigma }(s,t)\left\vert _{s=0}\right. .
\label{eq2.2}
\end{equation}

Multiplying Eq. (\ref{eq1.1}) by $s^{n}$ and summing over the $n$ yield

\begin{equation}
G_{\mu ,\nu }^{\gamma ,\sigma }(s,t)=\Gamma (\nu )E_{\mu ,\nu }^{\gamma
}(\lambda _{\sigma }t^{\sigma }(s-1)).  \label{eq2.3}
\end{equation}

Replacing $\lambda _{\sigma }t^{\sigma }$ with $x^{\sigma }$ gives us the
probability generating function for the distribution defined by Eq.(\ref%
{eq1.6a}),

\begin{equation}
g_{\mu ,\nu }^{\gamma ,\sigma }(s,x)=\Gamma (\nu )E_{\mu ,\nu }^{\gamma
}(x^{\sigma }(s-1)).  \label{eq2.3a}
\end{equation}

For $\mu =\nu =\gamma =1$ we obtain an expression for probability generating
function $G_{1,1}^{1,\sigma }(s,t)$ of the process introduced by Eq.(\ref%
{eq1.10_a})

\begin{equation}
G_{1,1}^{1,\sigma }(s,t)=\text{exp}(\lambda _{\sigma }t^{\sigma }(s-1)),
\label{eq2.4}
\end{equation}

here we have taken into account that at $\mu =\nu =\gamma =1$ the
three-parameter Mittag-Leffler function $E_{\mu ,\nu }^{\gamma }(z)$ becomes
the exponential function exp$(z)$.

For $\mu =\nu =\gamma =\sigma =1$ we obtain the well-known expression for
the probability generating function $G_{1,1}^{1,1}(s,t)$ of the Poisson
process

\begin{equation}
G_{1,1}^{1,1}(s,t)=\text{exp}(\lambda _{1}t(s-1)).  \label{eq2.4a}
\end{equation}

For $0<\mu \leq 1$, $\nu =\gamma =1$ and $\sigma =\mu $ we obtain an
expression for the probability generating function $G_{\mu ,1}^{1,\mu }(s,t)$
of the fractional Poisson process \cite{Laskin1}

\begin{equation}
G_{\mu ,1}^{1,\mu }(s,t)=E_{\mu }(\lambda _{\mu }t^{\mu }(s-1)),
\label{eq2.5}
\end{equation}

here we took into account that for $0<\mu \leq 1$, $\nu =\gamma =1$ the
three-parameter Mittag-Leffler function $E_{\mu ,1}^{1}(z)$ becomes the
Mittag-Leffler function $E_{\mu }(z)$ defined by its power series \cite{ML}, 
\cite{Erdelyi}

\begin{equation}
E_{\mu }(z)=\sum\limits_{m=0}^{\infty }\frac{z^{m}}{\Gamma (\mu m+1)}.
\label{eq2.5a}
\end{equation}

\subsection{Moment generating function}

An equation for the moment of any integer order of the probability
distribution function introduced by Eq. (\ref{eq1.1}) can be easily found
using the moment generating function $M_{\mu ,\nu }^{\gamma ,\sigma }(s,t)$,

\begin{equation}
M_{\mu ,\nu }^{\gamma ,\sigma }(s,t)=\sum\limits_{n=0}^{\infty
}e^{-sn}P_{\mu ,\nu }^{\gamma ,\sigma }(n,t).  \label{eq2.6}
\end{equation}

Thus, for the moment of $k^{\text{th}}$ order we have

\begin{equation}
<(n_{\mu ,\nu }^{\gamma ,\sigma })^{k}>=(-1)^{k}\frac{\partial ^{k}}{%
\partial s^{k}}M_{\mu ,\nu }^{\gamma ,\sigma }(s,t)\left\vert _{s=0}\right.
=\sum\limits_{n=0}^{\infty }n^{k}P_{\mu ,\nu }^{\gamma ,\sigma }(n,t).
\label{eq2.7}
\end{equation}

Multiplying Eq. (\ref{eq1.1}) by $e^{-sn}$ and summing over the $n$ we get

\begin{equation}
M_{\mu ,\nu }^{\gamma ,\sigma }(s,t)=\Gamma (\nu )E_{\mu ,\nu }^{\gamma
}(\lambda _{\sigma }t^{\sigma }(e^{-s}-1)).  \label{eq2.8a}
\end{equation}

Replacing $\lambda _{\sigma }t^{\sigma }$ with $x^{\sigma }$ gives us the
moment generating function for the distribution defined by Eq.(\ref{eq1.6a}),

\begin{equation}
m_{\mu ,\nu }^{\gamma ,\sigma }(s,x)=\Gamma (\nu )E_{\mu ,\nu }^{\gamma
}(x^{\sigma }(e^{-s}-1)).  \label{eq2.8b}
\end{equation}

For $\mu =\nu =\gamma =1$ we obtain an expression for the moment generating
function $M_{1,1}^{1,\sigma }(s,t)$ of the process introduced by Eq.(\ref%
{eq1.10_a})

\begin{equation*}
M_{1,1}^{1,\sigma }(s,t)=\exp (\lambda _{\sigma }t^{\sigma }(e^{-s}-1)).
\end{equation*}

For $\mu =\nu =\gamma =\sigma =1$ we obtain the well-known expression for
the moment generating function $M_{1,1}^{1,1}(s,t)$ of the Poisson process

\begin{equation*}
M_{1,1}^{1,1}(s,t)=\exp (\lambda _{1}t(e^{-s}-1)).
\end{equation*}

For $0<\mu \leq 1$, $\nu =\gamma =1$ and $\sigma =\mu $ we obtain an
expression for moment generating function $M_{\mu ,1}^{1,\mu }(s,t)$ of the
fractional Poisson process \cite{Laskin1}%
\begin{equation*}
M_{\mu ,1}^{1,\mu }(s,t)=E_{\mu }(\lambda _{\mu }t^{\mu }(e^{-s}-1)).
\end{equation*}

\section{Mean and variance}

The mean $<n_{\mu ,\nu }^{\gamma ,\sigma }>$ of the process defined by Eq.(%
\ref{eq1.1}) can be calculated directly using Eq.(\ref{eq1.14}),

\begin{equation}
<n_{\mu ,\nu }^{\gamma ,\sigma }>=\sum\limits_{n=0}^{\infty }nP_{\mu ,\nu
}^{\gamma ,\sigma }(n,t)=\gamma \Gamma (\nu )\frac{\lambda _{\sigma
}t^{\sigma }}{\Gamma (\mu +\nu )}.  \label{eq2.10}
\end{equation}

Alternatively, the mean can be calculated using Eqs.(\ref{eq2.7}) and (\ref%
{eq2.8a}). Indeed, we have the following chain of transformations

\begin{equation*}
<n_{\mu ,\nu }^{\gamma ,\sigma }>=(-1)\frac{\partial }{\partial s}M_{\mu
,\nu }^{\gamma ,\sigma }(s,t)\left\vert _{s=0}\right. =
\end{equation*}%
\begin{equation}
(-1)\Gamma (\nu )\frac{\partial }{\partial s}E_{\mu ,\nu }^{\gamma }(\lambda
_{\sigma }t^{\sigma }(e^{-s}-1))\left\vert _{s=0}\right. =  \label{eq2.11}
\end{equation}%
\begin{equation*}
(-1)\Gamma (\nu )\frac{\partial }{\partial s}\sum\limits_{m=0}^{\infty }%
\frac{(\gamma )_{m}}{m!}\frac{(\lambda _{\sigma }t^{\sigma
})^{m}(e^{-s}-1)^{m}}{\Gamma (\mu m+\nu )}\left\vert _{s=0}\right. =\gamma
\Gamma (\nu )\frac{\lambda _{\sigma }t^{\sigma }}{\Gamma (\mu +\nu )},
\end{equation*}

where we used the fact that only term $m=1$ contributes to the sum over $m$
after differentiating with respect to $s$ and then setting $s=0$.

The moment of second order $<(n_{\mu ,\nu }^{\gamma ,\sigma })^{2}>$ can be
calculated as follows

\begin{equation}
<(n_{\mu ,\nu }^{\gamma ,\sigma })^{2}>=(-1)^{2}\frac{\partial ^{2}}{%
\partial s^{2}}M_{\mu ,\nu }^{\gamma ,\sigma }(s,t)\left\vert _{s=0}\right. =
\label{eq2.12}
\end{equation}

\begin{equation*}
\gamma \Gamma (\nu )\frac{\lambda _{\sigma }t^{\sigma }}{\Gamma (\mu +\nu )}%
+(\gamma )_{2}\Gamma (\nu )\frac{(\lambda _{\sigma }t^{\sigma })^{2}}{\Gamma
(2\mu +\nu )}.
\end{equation*}

It is easy to see that $<(n_{\mu ,\nu }^{\gamma ,\sigma })^{2}>$ can be
represented in the form

\begin{equation}
<(n_{\mu ,\nu }^{\gamma ,\sigma })^{2}>=<n_{\mu ,\nu }^{\gamma ,\sigma
}>+<(n_{\mu ,\nu }^{\gamma ,\sigma })>^{2}(1+\frac{1}{\gamma })\frac{\text{B}%
(\mu +\nu ,\mu +\nu )}{\text{B}(2\mu +\nu ,\nu )},  \label{eq2.13}
\end{equation}

where B$(\alpha ,\beta )$ is the Beta-function defined in terms of the Gamma
function as%
\begin{equation}
\text{B}(\alpha ,\beta )=\frac{\Gamma (\alpha )\Gamma (\beta )}{\Gamma
(\alpha +\beta )}.  \label{eq2.14}
\end{equation}

Then variance $Var_{\mu ,\nu }^{\gamma ,\sigma }$ for the probability
distribution function $P_{\mu ,\nu }^{\gamma ,\sigma }(n,t)$ is

\begin{equation}
Var_{\mu ,\nu }^{\gamma ,\sigma }=<(n_{\mu ,\nu }^{\gamma ,\sigma
})^{2}>-<n_{\mu ,\nu }^{\gamma ,\sigma }>^{2}=  \label{eq2.15}
\end{equation}

\begin{equation*}
<n_{\mu ,\nu }^{\gamma ,\sigma }>+<n_{\mu ,\nu }^{\gamma ,\sigma
}>^{2}\left\{ (1+\frac{1}{\gamma })\frac{\text{B}(\mu +\nu ,\mu +\nu )}{%
\text{B}(2\mu +\nu ,\nu )}-1\right\} .
\end{equation*}

Note that when $\mu =\nu =\gamma =\sigma =1$ Eqs.(\ref{eq2.10}) - (\ref%
{eq2.14})\ are transformed into the well-known equations for the Poisson
probability distribution, and at $0<\mu \leq 1$, $\nu =\gamma =1$ and $%
\sigma =\mu $ Eqs.(\ref{eq2.10}) - (\ref{eq2.14})\ are transformed into the
equations for fractional Poisson probability distribution \cite{Laskin1}.

For $\mu =\nu =\gamma =1$ we obtain new equations for the mean $<n_{\sigma }>
$ and the second order moment $<(n_{\sigma })^{2}>$ of the probability
distribution given by Eq.(\ref{eq1.10_a})

\begin{equation}
<n_{\sigma }>=<n_{1,1}^{1,\sigma }>=\sum\limits_{n=0}^{\infty
}nP_{1,1}^{1,\sigma }(n,t)=\lambda _{\sigma }t^{\sigma },  \label{eq2.61}
\end{equation}

and

\begin{equation*}
<(n_{\sigma })^{2}>=<(n_{\mu ,\nu }^{\gamma ,\sigma
})^{2}>=\sum\limits_{n=0}^{\infty }n^{2}P_{1,1}^{1,\sigma }(n,t)=\lambda
_{\sigma }t^{\sigma }+(\lambda _{\sigma }t^{\sigma })^{2}.
\end{equation*}

Then the variance $Var_{\sigma }$ of the probability distribution Eq.(\ref%
{eq1.10_a}) has the form

\begin{equation*}
Var_{\sigma }\equiv Var_{1,1}^{1,\sigma }=<(n_{\mu ,\nu }^{\gamma ,\sigma
})^{2}>-<n_{\mu ,\nu }^{\gamma ,\sigma }>^{2}=\lambda _{\sigma }t^{\sigma }.
\end{equation*}

\section{Waiting time distribution of the fractional generalized counting
process}

We introduce waiting time probability distribution function $\psi _{\mu ,\nu
}^{\gamma ,\sigma }(\tau )$ of a new fractional generalized counting process 
$P_{\mu ,\nu }^{\gamma ,\sigma }(n,t)$ in the following way

\begin{equation}
\psi _{\mu ,\nu }^{\gamma ,\sigma }(\tau )=-\frac{d}{d\tau }P_{\mu ,\nu
}^{\gamma ,\sigma }(\tau ),  \label{eq3.1}
\end{equation}

where $P_{\mu ,\nu }^{\gamma ,\sigma }(\tau )$ is the probability that a
given waiting time is greater or equal to $\tau $

\begin{equation}
P_{\mu ,\nu }^{\gamma ,\sigma }(\tau )=1-\sum\limits_{n=1}^{\infty }P_{\mu
,\nu }^{\gamma ,\sigma }(n,\tau )=P_{\mu ,\nu }^{\gamma ,\sigma }(0,\tau
)=\Gamma (\nu )E_{\mu ,\nu }^{\gamma }(-\lambda _{\sigma }^{\sigma }\tau
^{\sigma }).  \label{eq3.2}
\end{equation}

Using Eqs.(\ref{eq3.1}) and (\ref{eq3.2}) we obtain the waiting time
probability distribution

\begin{equation}
\psi _{\mu ,\nu }^{\gamma ,\sigma }(\tau )=-\Gamma (\nu )\frac{d}{d\tau }%
E_{\mu ,\nu }^{\gamma }(-\lambda _{\sigma }\tau ^{\sigma })=\gamma \sigma
\lambda _{\sigma }\tau ^{\sigma -1}E_{\mu ,\nu +\mu }^{\gamma +1}(-\lambda
_{\sigma }\tau ^{\sigma }),  \label{eq3.3}
\end{equation}

where Eq.(\ref{eq1.11}) is taken into account.

For $\mu =\nu =\gamma =\sigma =1$ the probability distribution function $%
\psi _{\mu ,\nu }^{\gamma ,\sigma }(\tau )$ becomes the well-known
exponential distribution for the Poisson process,

\begin{equation}
\psi (\tau )\equiv \psi _{1,1}^{1,1}(\tau )=\lambda _{1}e^{-\lambda _{1}\tau
}.  \label{eq3.4}
\end{equation}

For $0<\mu \leq 1$, $\nu =\gamma =1$ and $\sigma =\mu $ we obtain the
waiting time distribution for the fractional Poisson process \cite{Laskin1}

\begin{equation}
\psi _{\mu }(\tau )\equiv \psi _{\mu ,1}^{1,\mu }(\tau )=\lambda _{\mu }\tau
^{\mu -1}E_{\mu ,\mu }(-\lambda _{\mu }\tau ^{\mu }).  \label{eq3.5}
\end{equation}

In the case when $\mu =$ $\nu =\gamma =1$ and $0<\sigma \leq 1$ we have

\begin{equation}
\psi _{\sigma }(\tau )=\psi _{1,1}^{1,\sigma }(\tau )=\sigma \lambda
_{\sigma }\tau ^{\sigma -1}\exp (-\lambda _{\sigma }\tau ^{\sigma }),\qquad
\tau >0,  \label{eq3.6}
\end{equation}

which is the waiting time probability distribution for the counting process
introduced by Eq.(\ref{eq1.10_a}).

One can see that $\psi _{\sigma }(\tau )$ is the well-known Weibull
distribution with the shape parameter $\sigma $ and the scale parameter $%
\lambda _{\sigma }$. Probability distribution $\psi _{\sigma }(\tau )$ was
used to fit "fat-tail" economic and natural data in \cite{LAHERRERE} and is
called the \textit{stretched exponential distribution}.

\subsection{Moments of a stretched exponential distribution}

The moments of the stretched exponential distribution $\psi _{\sigma }(\tau
) $ are easy to calculate. For the mean $<\tau _{\sigma }>$ we have

\begin{equation}
<\tau _{\sigma }>=\int\limits_{0}^{\infty }d\tau \tau \psi _{\sigma }(\tau )=%
\frac{1}{\lambda _{\sigma }^{1/\sigma }}\int\limits_{0}^{\infty
}dyy^{1/\sigma }e^{-y}=\frac{1}{\sigma \lambda _{\sigma }^{1/\sigma }}\Gamma
(1/\sigma ),  \label{eq3.7}
\end{equation}

and for the moment of the second order $<\tau _{\sigma }^{2}>,$

\begin{equation}
<\tau _{\sigma }^{2}>=\int\limits_{0}^{\infty }d\tau \tau ^{2}\psi _{\sigma
}(\tau )=\frac{2}{\sigma \lambda _{\sigma }^{2/\sigma }}\Gamma (2/\sigma ).
\label{eq3.8}
\end{equation}

Then the variance of the stretched exponential distribution is

\begin{equation}
Var_{\sigma }=<\tau _{\sigma }^{2}>-<\tau _{\sigma }>^{2}=  \label{eq3.9}
\end{equation}

\begin{equation*}
\frac{1}{\sigma \lambda _{\sigma }^{2/\sigma }}\left\{ 2\Gamma (2/\sigma )-%
\frac{1}{\sigma }\Gamma ^{2}(1/\sigma )\right\} .
\end{equation*}

Using the Legendre duplication formula for the Gamma function

\begin{equation*}
\Gamma (2z)=\frac{1}{\sqrt{\pi }}2^{2z-1}\Gamma (z)\Gamma (z+\frac{1}{2}),
\end{equation*}

we express the variance of the stretched exponential distribution in an
alternative form

\begin{equation}
Var_{\sigma }=\frac{\Gamma ^{2}(1/\sigma )}{\sigma ^{2}\lambda _{\sigma
}^{2/\sigma }}\left\{ \frac{(2^{2/\sigma }\sigma )^{2}}{\text{B}(1/\sigma
,1/2)}-1\right\} ,  \label{eq3.10}
\end{equation}

where B$(1/\sigma ,1/2)$ is the Beta function defined by Eq.(\ref{eq2.14}).

It is easy to see that the formula for the moment of any order $k$ has the
form

\begin{equation*}
<\tau _{\sigma }^{k}>=\int\limits_{0}^{\infty }d\tau \tau ^{k}\psi _{\sigma
}(\tau )=\frac{k}{\sigma }\frac{\Gamma (k/\sigma )}{\lambda _{\sigma
}^{k/\sigma }}.
\end{equation*}

\section{Applications of fractional probability distributions}

\subsection{Fractional generalized compound process}

We call stochastic process $\{X_{\mu ,\nu }^{\gamma ,\sigma }(t)$, $t\geq
0\} $ a fractional generalized compound process if it is represented by

\begin{equation}
X_{\mu ,\nu }^{\gamma ,\sigma }(t)=\sum\limits_{i=1}^{N_{\mu ,\nu }^{\gamma
,\sigma }(t)}Y_{i},  \label{eq4.1}
\end{equation}

where $\{N_{\mu ,\nu }^{\gamma ,\sigma }(t)$, $t\geq 0\}$ is a fractional
generalized counting process whose probability distribution function is
introduced by Eq.(\ref{eq1.1}), and $\{Y_{i}$, $i=1,2,...\}$ are independent
and identically distributed random variables with probability distribution
function $p(Y_{i})$ for each $Y_{i}$. It is assumed that the process $%
\{N_{\mu ,\nu }^{\gamma ,\sigma }(t)$, $t\geq 0\}$ and the sequence of
random variables $\{Y_{i}$, $i=1,2,...\}$ are independent.

The moment generating function $J(s,t)$ of fractional generalized compound
process is defined as follows

\begin{equation}
J_{\mu ,\nu }^{\gamma ,\sigma }(s,t)=<\exp \{sX_{\mu ,\nu }^{\gamma ,\sigma
}(t)\}>_{Y,N_{\mu ,\nu }^{\gamma ,\sigma }(t)},  \label{eq4.2}
\end{equation}

where $<...>_{Y_{i},N_{\mu ,\nu }^{\gamma ,\sigma }(t)}$ stands for average
which includes two statistically independent averaging procedures:

1. Averaging $<...>_{Y}$ over independent random variables $Y_{i}$, defined
by

\begin{equation}
<...>_{Y}=\int
dY_{1}...dY_{n}p(Y_{1})...p(Y_{n})...=\prod\limits_{i=1}^{n}\int
dY_{i}p(Y_{i})....  \label{eq4.3}
\end{equation}

2. Averaging over random counts defined by

\begin{equation}
<...>_{N_{\mu ,\nu }^{\gamma ,\sigma }(t)}=\sum\limits_{n=0}^{\infty }P_{\mu
,\nu }^{\gamma ,\sigma }(n,t)...,  \label{eq4.4}
\end{equation}

where $P_{\mu ,\nu }^{\gamma ,\sigma }(n,t)$ is given by Eq.(\ref{eq1.1}).

One can see from Eq.(\ref{eq4.2}) that $k^{\mathrm{th}}$ order moment of
fractional generalized compound process $X_{\mu ,\nu }^{\gamma ,\sigma }(t)$
is obtained by differentiating $J_{\mu }(s,t)$ $k$ times with respect to $s$
and then setting $s=0$, that is

\begin{equation}
<(X_{\mu ,\nu }^{\gamma ,\sigma }(t))^{k}>_{Y,N_{\mu ,\nu }^{\gamma ,\sigma
}(t)}=\frac{\partial ^{k}}{\partial s^{k}}J_{\mu ,\nu }^{\gamma ,\sigma
}(s,t)_{|s=0}.  \label{eq4.5}
\end{equation}

To obtain equation for the moment generating function $J_{\mu ,\nu }^{\gamma
,\sigma }(s,t)$ we apply Eqs.(\ref{eq4.3}) and (\ref{eq4.4}) to Eq.(\ref%
{eq4.2})

\begin{equation}
J_{\mu ,\nu }^{\gamma ,\sigma }(s,t)=\sum\limits_{n=0}^{\infty }<\exp
\{sX_{\mu ,\nu }^{\gamma ,\sigma }(t)|N_{\mu ,\nu }^{\gamma ,\sigma
}(t)=n\}>_{Y}P_{\mu ,\nu }^{\gamma ,\sigma }(n,t)=  \label{eq4.6}
\end{equation}

\begin{equation*}
\sum\limits_{n=0}^{\infty }(\prod\limits_{i=1}^{n}\int dY_{i}p(Y_{i})(\exp
\{sY_{i}|N_{\mu ,\nu }^{\gamma ,\sigma }(t)=n\})P_{\mu ,\nu }^{\gamma
,\sigma }(n,t)=
\end{equation*}

\begin{equation*}
\sum\limits_{n=0}^{\infty }<\exp \{s(Y_{1}+...+Y_{n})\}>_{Y}\times \Gamma
(\nu )\frac{(-z)^{n}}{n!}\frac{d^{n}}{dz^{n}}\left( E_{\mu ,\nu }^{\gamma
}(z)\right) |_{z=-\lambda _{\sigma }t^{\sigma }}=
\end{equation*}

\begin{equation*}
\Gamma (\nu )\sum\limits_{n=0}^{\infty }\prod\limits_{i=1}^{n}\int
dY_{i}p(Y_{i})\exp \{sY_{i}\})\times \frac{(-z)^{n}}{n!}\frac{d^{n}}{dz^{n}}%
\left( E_{\mu ,\nu }^{\gamma }(z)\right) |_{z=-\lambda _{\sigma }t^{\sigma
}},
\end{equation*}

where we used the independence of $\{Y_{1},Y_{2},...\}$ and $\{N_{\mu ,\nu
}^{\gamma ,\sigma }(t)$, $t\geq 0\}$ and the independence of the $Y_{i}$'s
between themselves. Hence, setting

\begin{equation}
g(s)=\int dY_{i}p(Y_{i})\exp \{sY_{i}\},  \label{eq4.7}
\end{equation}

for the moment generating function of random variables $Y_{i}$, we obtain
from Eq.(\ref{eq4.6}) the moment generating function $J_{\mu ,\nu }^{\gamma
,\sigma }(s,t)$ of the fractional generalized compound process

\begin{equation}
J_{\mu ,\nu }^{\gamma ,\sigma }(s,t)=\Gamma (\nu )\sum\limits_{n=0}^{\infty
}g^{n}(s)\times \frac{(-z)^{n}}{n!}\frac{d^{n}}{dz^{n}}\left( E_{\mu ,\nu
}^{\gamma }(z)\right) |_{z=-\lambda _{\sigma }t^{\sigma }}=  \label{eq.4.8}
\end{equation}%
\begin{equation*}
\Gamma (\nu )E_{\mu ,\nu }^{\gamma }(\lambda _{\sigma }t^{\sigma }(g(s)-1)).
\end{equation*}

By differentiating the above with respect to $s$ and setting $s=0$, it is
easy to obtain the mean of the fractional generalized compound process,

\begin{equation}
<X_{\mu ,\nu }^{\gamma ,\sigma }(t)>_{Y_{i},N_{\mu ,\nu }^{\gamma ,\sigma
}(t)}=\frac{\partial }{\partial s}J_{\mu ,\nu }^{\gamma ,\sigma
}(s,t)_{|s=0}=<Y>_{Y}\gamma \Gamma (\nu )\frac{\lambda _{\sigma }t^{\sigma }%
}{\Gamma (\mu +\nu )},  \label{eq4.9}
\end{equation}

which is a manifestation of independency between fractional counting process 
$\{N_{\mu ,\nu }^{\gamma ,\sigma }(t)$, $t\geq 0\}$ and random variables $%
Y_{i}.$

\subsection{Fractional generalized Bell polynomials and fractional
generalized Bell numbers}

Based on the new probability distribution function defined by Eq.(\ref{eq1.1}%
) we introduce a new fractional generalization of the Bell polynomials

\begin{equation}
B_{\mu ,\nu }^{\gamma }(x,m)=\Gamma (\nu )\sum\limits_{n=0}^{\infty }n^{m}%
\frac{x^{n}}{n!}E_{\mu ,\nu }^{\gamma ^{(n)}}(-x),\qquad B_{\mu ,\nu
}^{\gamma }(x,0)=1,  \label{eq5.1a}
\end{equation}

where the notation has been introduced

\begin{equation}
E_{\mu ,\nu }^{\gamma ^{(n)}}(-x)=\frac{d^{n}}{dz^{n}}\left( E_{\mu ,\nu
}^{\gamma }(z)\right) |_{z=-x},  \label{eq5.2}
\end{equation}

and the parameters $\mu $, $\nu $\ and $\gamma $ satisfy the conditions
given by Eq.(\ref{eq1.4}). We call $B_{\mu ,\nu }^{\gamma }(x,m)$ as the
fractional generalized Bell polynomials of $m^{\mathrm{th}}$ order.

The polynomials $B_{\mu ,\nu }^{\gamma }(x,m)$ are related to the well-known
Bell polynomials $B(x,m)$ \cite{Bell} by

\begin{equation}
B(x,m)=B_{\mu ,\nu }^{\gamma }(x,m)\left\vert _{\mu =\nu =\gamma =1}\right.
=B_{1,1}^{1}(x,m)=e^{\text{-}x}\sum\limits_{n=0}^{\infty }n^{m}\frac{x^{n}}{%
n!}.  \label{eq5.3}
\end{equation}

The polynomials $B_{\mu ,\nu }^{\gamma }(x,m)$ are related to the fractional
Bell polynomials $B_{\mu }(x,m)$ \cite{Laskin2} by

\begin{equation}
B_{\mu }(x,m)=B_{\mu ,\nu }^{\gamma }(x,m)\left\vert _{\nu =\gamma
=1}\right. =B_{\mu ,1}^{1}(x,m)=\sum\limits_{n=0}^{\infty }n^{m}\frac{x^{n}}{%
n!}E_{\mu ,1}^{1^{(n)}}(-x).  \label{eq5.4}
\end{equation}

From Eq.(\ref{eq5.1a}) at $x=1$ we come to new numbers, which we call
fractional generalized Bell numbers

\begin{equation}
B_{\mu ,\nu }^{\gamma }(m)=B_{\mu ,\nu }^{\gamma }(x,m)|_{x=1}=\Gamma (\nu
)\sum\limits_{n=0}^{\infty }n^{m}\frac{x^{n}}{n!}E_{\mu ,\nu }^{\gamma
^{(n)}}(-1).  \label{eq5.5}
\end{equation}

The new formula Eq.(\ref{eq5.5}) is a generalization of the so-called Dobi%
\'{n}ski relation known since 1877 \cite{Dobinski} for the Bell numbers%
\footnote{%
The Bell number is the number of ways in which a set of $m$ elements can be
partitioned into non-empty subsets.} $B(m)$

\begin{equation}
B(m)=B_{1,1}^{1}(m)=e^{\text{-}1}\sum\limits_{n=0}^{\infty }\frac{n^{m}}{n!}.
\label{eq5.6}
\end{equation}

When $\nu =\gamma =1$, we recover the equation for fractional Bell numbers
originally found in \cite{Laskin2}

\begin{equation}
B_{\mu }(m)=B_{\mu ,1}^{1}(m)=\sum\limits_{n=0}^{\infty }\frac{n^{m}}{n!}%
E_{\mu }^{^{(n)}}(-1).  \label{eq5.6a}
\end{equation}

As an example, we present the first two newly introduced fractional
generalized Bell polynomials

\begin{equation*}
B_{\mu ,\nu }^{\gamma }(x,1)=\frac{\gamma \Gamma (\nu )}{\Gamma (\mu +\nu )}%
x,\qquad B_{\mu ,\nu }^{\gamma }(x,2)=\frac{(\gamma )_{2}\Gamma (\nu )}{%
\Gamma (2\mu +\nu )}x^{2}+\frac{\gamma \Gamma (\nu )}{\Gamma (\mu +\nu )}x,
\end{equation*}

and the first two fractional generalized Bell numbers

\begin{equation*}
B_{\mu ,\nu }^{\gamma }(1)=\frac{\gamma \Gamma (\nu )}{\Gamma (\mu +\nu )}%
,\qquad B_{\mu ,\nu }^{\gamma }(2)=\frac{(\gamma )_{2}\Gamma (\nu )}{\Gamma
(2\mu +\nu )}+\frac{\gamma \Gamma (\nu )}{\Gamma (\mu +\nu )}.
\end{equation*}

The alternative representations for fractional generalized Bell polynomials $%
B_{\mu ,\nu }^{\gamma }(x,m)$ and fractional generalized Bell numbers $%
B_{\mu ,\nu }^{\gamma }(m)$ are%
\begin{equation}
B_{\mu ,\nu }^{\gamma }(x,m)=\Gamma (\nu )\sum\limits_{n=0}^{\infty }(\gamma
)_{n}n^{m}\frac{x^{n}}{n!}\sum\limits_{k=0}^{\infty }\frac{(\gamma +n)_{k}}{%
k!\Gamma (\mu (n+k)+\nu )}(-x)^{k},  \label{eq5.7}
\end{equation}

\begin{equation}
B_{\mu ,\nu }^{\gamma }(m)=B_{\mu ,\nu }^{\gamma }(x,m)|_{x=1}=\Gamma (\nu
)\sum\limits_{n=0}^{\infty }(\gamma )_{n}\frac{n^{m}}{n!}\sum\limits_{k=0}^{%
\infty }\frac{(-1)^{k}(\gamma +n)_{k}}{k!\Gamma (\mu (n+k)+\nu )},
\label{eq5.8}
\end{equation}

where we used Eq.(\ref{eq1.14}).

\subsubsection{Generating function of fractional generalized Bell polynomials%
}

The generating function of the fractional generalized Bell polynomials is

\begin{equation}
F_{\mu ,\nu }^{\gamma }(s,x)=\sum\limits_{m=0}^{\infty }\frac{s^{m}}{m!}%
B_{\mu ,\nu }^{\gamma }(x,m).  \label{eq6.1}
\end{equation}

To get the polynomial $B_{\mu ,\nu }^{\gamma }(x,m)$ we differentiate $%
F_{\mu ,\nu }^{\gamma }(s,x)$ $m$ times with respect to $s$, and then let $%
s=0$. That is,

\begin{equation}
B_{\mu ,\nu }^{\gamma }(x,m)=\frac{\partial ^{m}}{\partial s^{m}}F_{\mu ,\nu
}^{\gamma }(s,x)\left\vert _{s=0}\right. .  \label{eq6.2}
\end{equation}

To find an explicit equation for $F_{\mu ,\nu }^{\gamma }(s,x)$, let's
substitute Eq.(\ref{eq5.1a}) into Eq.(\ref{eq6.1}) and evaluate the sum over 
$m$. As a result we have

\begin{equation}
F_{\mu ,\nu }^{\gamma }(s,x)=\Gamma (\nu )E_{\mu ,\nu }^{\gamma
}(x(e^{s}-1)).  \label{eq6.3}
\end{equation}

If we put $x=1$ in Eq.(\ref{eq6.3}), then we immediately come to the
generating function $\mathcal{F}_{\mu ,\nu }^{\gamma }(s)$ of the fractional
generalized Bell numbers $B_{\mu ,\nu }^{\gamma }(m)$ introduced by Eq.(\ref%
{eq5.5})

\begin{equation}
\mathcal{F}_{\mu ,\nu }^{\gamma }(s)=F_{\mu ,\nu }^{\gamma
}(s,x)|_{x=1}=\Gamma (\nu )E_{\mu ,\nu }^{\gamma }(e^{s}-1).  \label{eq6.4}
\end{equation}

To get the fractional generalized Bell numbers $B_{\mu ,\nu }^{\gamma }(m)$
we differentiate $\mathcal{F}_{\mu ,\nu }^{\gamma }(s)$ $m$ times with
respect to $s$, and then set $s=0$. That is,

\begin{equation}
B_{\mu ,\nu }^{\gamma }(m)=\frac{\partial ^{m}}{\partial s^{m}}\mathcal{F}%
_{\mu ,\nu }^{\gamma }(s)\left\vert _{s=0}\right. .  \label{eq6.5}
\end{equation}

\subsection{Fractional generalized Stirling numbers of the second kind}

We introduce the fractional generalized Stirling numbers of the second kind $%
S_{\mu ,\nu }^{\gamma }(x,l)$ by means of the equation

\begin{equation}
B_{\mu ,\nu }^{\gamma }(x,m)=\sum\limits_{l=0}^{m}S_{\mu ,\nu }^{\gamma
}(m,l)x^{l},  \label{eq8.1}
\end{equation}

\begin{equation*}
S_{\mu ,\nu }^{\gamma }(m,0)=\delta _{m,0},\qquad S_{\mu ,\nu }^{\gamma
}(m,l)=0,\ l=m+1,\ m+2,\ ...,
\end{equation*}

where $B_{\mu ,\nu }^{\gamma }(x,m)$ is fractional generalized Bell
polynomials given by Eq.(\ref{eq5.1a}) and the parameters $\mu $, $\nu $\
and $\gamma $ satisfy the conditions (\ref{eq1.4}). For $\mu =\nu =\gamma =1$
the fractional generalized Stirling numbers of the second kind $S_{\mu ,\nu
}^{\gamma }(m,l)$ becomes the well-known Stirling numbers of the second kind%
\footnote{%
The number of ways to partition a set of $m$ elements into $l$ non-empty
sets.} $S(m,l)=S_{1,1}^{1}(m,l)$, see \cite{Stirling}.

At $x=1$, when the fractional generalized Bell polynomials $B_{\mu ,\nu
}^{\gamma }(x,m)$ become the fractional generalized Bell numbers, $B_{\mu
,\nu }^{\gamma }(m)=B_{\mu ,\nu }^{\gamma }(x,m)\left\vert _{x=1}\right. $,
Eq.(\ref{eq8.1}) goes into the expression of the fractional generalized Bell
numbers in terms of fractional generalized Stirling numbers of the second
kind

\begin{equation}
B_{\mu ,\nu }^{\gamma }(m)=\sum\limits_{l=0}^{m}S_{\mu ,\nu }^{\gamma }(m,l).
\label{eq8.2}
\end{equation}

To find a generating function of the fractional generalized Stirling numbers 
$S_{\mu ,\nu }^{\gamma }(m,l)$, let's expand the generating function $F_{\mu
,\nu }^{\gamma }(s,x)$ given by Eq.(\ref{eq6.1}). Upon substituting $B_{\mu
}(x,m)$ from Eq.(\ref{eq8.1}) we have the following chain of transformations

\begin{equation*}
F_{\mu ,\nu }^{\gamma }(s,x)=\sum\limits_{m=0}^{\infty }\frac{s^{m}}{m!}%
\left( \sum\limits_{l=0}^{m}S_{\mu ,\nu }^{\gamma }(m,l)x^{l}\right) =
\end{equation*}

\begin{equation}
\sum\limits_{m=0}^{\infty }\frac{s^{m}}{m!}\left( \sum\limits_{l=0}^{\infty
}\theta (m-l)S_{\mu ,\nu }^{\gamma }(m,l)x^{l}\right)
=\sum\limits_{l=0}^{\infty }\left( \sum\limits_{m=l}^{\infty }S_{\mu ,\nu
}^{\gamma }(m,l)\frac{s^{m}}{m!}\right) x^{l},  \label{eq8.3}
\end{equation}

here $\theta (l)$ is the Heaviside step function,

\begin{equation}
\theta (l)=\mid \QATOP{1,\quad \ \mathrm{if\ }l\geq 0}{0,\ \quad \mathrm{if\ 
}l<0}.  \label{eq8.4}
\end{equation}

On the other hand, from Eq.(\ref{eq6.3}), for $F_{\mu ,\nu }^{\gamma }(s,x)$
we have

\begin{equation}
F_{\mu ,\nu }^{\gamma }(s,x)=\Gamma (\nu )E_{\mu ,\nu }^{\gamma
}(x(e^{s}-1))=\Gamma (\nu )\sum\limits_{l=0}^{\infty }\frac{(\gamma
)_{l}(x(e^{s}-1))^{l}}{l!\Gamma (\mu l+\nu )}.  \label{eq8.4a}
\end{equation}

Comparing Eq.(\ref{eq8.4a}) and Eq.(\ref{eq8.3}), we come to the conclusion
that

\begin{equation}
\sum\limits_{m=l}^{\infty }S_{\mu ,\nu }^{\gamma }(m,l)\frac{s^{m}}{m!}%
=\Gamma (\nu )\frac{(\gamma )_{l}(e^{s}-1)^{l}}{l!\Gamma (\mu l+\nu )}%
,\qquad l=0,1,2,....  \label{eq8.5}
\end{equation}

Let us introduce two generating functions $\mathcal{G}_{\mu ,\nu }^{\gamma
}(s,l)$ and $\mathcal{F}_{\mu ,\nu }^{\gamma }(s,t)$ of the fractional
generalized Stirling numbers of the second kind,

\begin{equation}
\mathcal{G}_{\mu ,\nu }^{\gamma }(s,l)=\sum\limits_{m=l}^{\infty }S_{\mu
,\nu }^{\gamma }(m,l)\frac{s^{m}}{m!}=\Gamma (\nu )\frac{(\gamma
)_{l}(e^{s}-1)^{l}}{l!\Gamma (\mu l+\nu )},  \label{eq8.6}
\end{equation}

and

\begin{equation}
\mathcal{F}_{\mu ,\nu }^{\gamma
}(s,t)=\sum\limits_{l=0}^{m}\sum\limits_{m=l}^{\infty }S_{\mu ,\nu }^{\gamma
}(m,l)\frac{s^{m}t^{l}}{m!}=\Gamma (\nu )\sum\limits_{l=0}^{\infty }\frac{%
(\gamma )_{l}(e^{s}-1)^{l}}{l!\Gamma (\mu l+\nu )}t^{l}=  \label{eq8.7}
\end{equation}

\begin{equation*}
\Gamma (\nu )E_{\mu ,\nu }^{\gamma }(t(e^{s}-1)).
\end{equation*}

Then for fractional generalized Stirling numbers of the second kind we obtain%
\begin{equation}
S_{\mu ,\nu }^{\gamma }(m,l)=\frac{\partial ^{m}}{\partial s^{m}}\mathcal{G}%
_{\mu ,\nu }^{\gamma }(s,l)\left\vert _{s=0}\right. =\Gamma (\nu )\frac{%
(\gamma )_{l}}{l!\Gamma (\mu l+\nu )}\frac{\partial ^{m}}{\partial s^{m}}%
(e^{s}-1)^{l}\left\vert _{s=0}\right. ,  \label{eq8.8}
\end{equation}

and

\begin{equation}
S_{\mu ,\nu }^{\gamma }(m,l)=\frac{1}{l!}\frac{\partial ^{m+l}}{\partial
s^{m}\partial t^{l}}\mathcal{F}_{\mu ,\nu }^{\gamma }(s,t)\left\vert
_{s=0},_{t=0}\right. =\Gamma (\nu )\frac{1}{l!}\frac{\partial ^{m+l}}{%
\partial s^{m}\partial t^{l}}E_{\mu ,\nu }^{\gamma
}(t(e^{s}-1))_{s=0},_{t=0},  \label{eq8.9}
\end{equation}

where $l\leq m$.

In the special case, for $\mu =\nu =\gamma =1$, equations (\ref{eq8.6}) and (%
\ref{eq8.7}) include well-know generating functions for standard Stirling
numbers of the second kind $S(m,l)$ (see, Eqs.(2.17) and (2.18) in \cite%
{Charalambides}),

\begin{equation}
\mathcal{G}_{1,1}^{1}(s,l)=\sum\limits_{m=l}^{\infty }S(m,l)\frac{s^{m}}{m!}=%
\frac{(e^{s}-1)^{l}}{l!},\qquad l=0,1,2,...,  \label{eq8.10}
\end{equation}

and

\begin{equation}
\mathcal{F}_{1,1}^{1}(s,t)=\sum\limits_{m=0}^{\infty
}\sum\limits_{l=0}^{m}S(m,l)\frac{s^{m}t^{l}}{m!}=\exp (t(e^{s}-1)).
\label{eq8.11}
\end{equation}

\subsubsection{Statistics of fractional generalized probability distribution}

We now use the fractional generalized Stirling numbers introduced by Eq.(\ref%
{eq8.1}) to obtain the moments of the fractional generalized process with
the probability distribution function $P_{\mu ,\nu }^{\gamma ,\sigma }(n,t)$
given by Eqs.(\ref{eq1.1}) and (\ref{eq1.13}). By the definition of the $m^{%
\mathrm{th}}$ order moment of the fractional generalized probability
distribution, we have

\begin{equation}
<(n_{\mu ,\nu }^{\gamma ,\sigma })^{m}>=\sum\limits_{n=0}^{\infty
}n^{m}P_{\mu ,\nu }^{\gamma ,\sigma }(n,t)=\Gamma (\nu
)\sum\limits_{n=0}^{\infty }n^{m}\frac{(-z)^{n}}{n!}\frac{d^{n}}{dz^{n}}%
\left( E_{\mu ,\nu }^{\gamma }(z)\right) |_{z=-\lambda _{\sigma }t^{\sigma
}},  \label{eq8.12}
\end{equation}

where the fractality parameters $\mu $, $\nu ,$\ $\gamma $ and $\sigma $
satisfy the conditions given by Eq.(\ref{eq1.4}). It is easy to see that $%
<(n_{\mu ,\nu }^{\gamma ,\sigma })^{m}>$ is the fractional generalized Bell
polynomial $B_{\mu ,\nu }^{\gamma }(\lambda _{\sigma }t^{\sigma },m)$
defined by Eq.(\ref{eq5.1a}). Then using the definition Eq.(\ref{eq8.1}) we
obtain

\begin{equation}
<(n_{\mu ,\nu }^{\gamma ,\sigma })^{m}>=\sum\limits_{l=0}^{m}S_{\mu ,\nu
}^{\gamma }(m,l)(\lambda _{\sigma }t^{\sigma })^{l},  \label{eq8.13}
\end{equation}

where $S_{\mu ,\nu }^{\gamma }(m,l)$ is a fractional generalized Stirling
number of the second kind introduced by Eq.(\ref{eq8.1}).

Therefore, the fractional generalized Stirling numbers $S_{\mu ,\nu
}^{\gamma }(m,l)$ of the second kind naturally appear in the power series
over $\lambda _{\sigma }t^{\sigma }$ for the $m^{\mathrm{th}}$ order moment
of the fractional generalized process introduced by Eq.(\ref{eq1.1}).

\subsection{New quantum coherent states}

To show how fractional generalized Bell polynomials appear in quantum
physics, we introduce a new family of quantum coherent states $|\varsigma
>_{\mu ,\nu }^{\gamma ,\sigma },$

\begin{equation}
|\varsigma >_{\mu ,\nu }^{\gamma ,\sigma }=\Gamma ^{1/2}(\nu
)\sum\limits_{n=0}^{\infty }\frac{\varsigma ^{\sigma n}}{\sqrt{n!}}(E_{\mu
,\nu }^{\gamma ^{(n)}}(-|\varsigma |^{2\sigma }))^{1/2}|n>,  \label{eq7.1}
\end{equation}

and the adjoint states $_{\mu ,\nu }^{\gamma ,\sigma }<\varsigma |$

\begin{equation}
_{\mu ,\nu }^{\gamma ,\sigma }<\varsigma |=\Gamma ^{1/2}(\nu
)\sum\limits_{n=0}^{\infty }<n|\frac{(\varsigma ^{\sigma \ast })^{n}}{\sqrt{%
n!}}(E_{\mu ,\nu }^{\gamma ^{(n)}}(-|\varsigma |^{2\sigma }))^{1/2},
\label{eq7.2}
\end{equation}

here 
\begin{equation}
E_{\mu ,\nu }^{\gamma ^{(n)}}(-|\varsigma |^{2\sigma })=\frac{d^{n}}{dz^{n}}%
E_{\mu ,\nu }^{\gamma ^{(n)}}(z)|_{z=-|\varsigma |^{2\sigma }},
\label{eq7.3}
\end{equation}

with $E_{\mu ,\nu }^{\gamma }(z)$ being the three-parameter Mittag-Leffler
function given by its power series Eq.(\ref{eq1.2}), complex number $%
\varsigma $ stands for labelling the new coherent states, the orthonormal
vectors $|n>$ are eigenvectors of the photon number operator $\overset{%
\wedge }{n}=a^{+}a$, $\overset{\wedge }{n}|n>=n|n>$, $<n|m>=\delta _{n,m}$
and the fractality parameters $\mu $, $\nu ,$\ $\gamma $ and $\sigma $
satisfy the conditions Eq.(\ref{eq1.4}). The operators $a^{+}$ and $a$ are
photon field creation and annihilation operators that satisfy the
Bose-Einstein commutation relation $[a,a^{+}]=aa^{+}-a^{+}a=\QTR{sl}{1}$,
and $[a^{+},a^{+}]=0$, $[a,a]=0$.

The action of the operators $a^{+}$ and $\ a$ on the number state $|n>$ reads

{}%
\begin{equation}
a^{+}|n>=\sqrt{n+1}|n+1>\qquad \text{and}\qquad a|n>=\sqrt{n}|n-1>.
\label{eq7.4}
\end{equation}

The motivation to introduce the coherent states $|\varsigma >_{\mu ,\nu
}^{\gamma ,\sigma }$ is the observation that the squared modulus $%
|<n|\varsigma >_{\mu ,\nu }^{\gamma ,\sigma }|^{2}$ of projection of the
coherent state $|\varsigma >_{\mu ,\nu }^{\gamma ,\sigma }$ onto the number
state $|n>$ gives us the probability $P_{\mu ,\nu }^{\gamma ,\sigma
}(n,\varsigma )$ that $n$ photons will be found in the coherent state $%
|\varsigma >_{\mu ,\nu }^{\gamma ,\sigma },$

\begin{equation}
P_{\mu ,\nu }^{\gamma ,\sigma }(n,\varsigma )=|<n|\varsigma >_{\mu ,\nu
}^{\gamma ,\sigma }|^{2}=\Gamma (\nu )\frac{|\varsigma |^{2\sigma n}}{n!}%
E_{\mu ,\nu }^{\gamma ^{(n)}}(-|\varsigma |^{2\sigma }),  \label{eq7.5}
\end{equation}

and

\begin{equation}
P_{\mu ,\nu }^{\gamma ,\sigma }(0,\varsigma )=\Gamma (\nu )E_{\mu ,\nu
}^{\gamma }(-|\varsigma |^{2\sigma }).  \label{eq7.6a}
\end{equation}

The alternative expression for $P_{\mu ,\nu }^{\gamma ,\sigma }(n,\varsigma
) $ is

\begin{equation}
P_{\mu ,\nu }^{\gamma ,\sigma }(n,\varsigma )=(\gamma )_{n}\Gamma (\nu )%
\frac{|\varsigma |^{2\sigma n}}{n!}\sum\limits_{k=0}^{\infty }\frac{(\gamma
+n)_{k}}{k!\Gamma (\mu (n+k)+\nu )}(-|\varsigma |^{2\sigma })^{k}.
\label{eq7.7a}
\end{equation}

It is easy to see, that the mean number of photons in the quantum state $%
|\varsigma >_{\mu ,\nu }^{\gamma ,\sigma }$ is

\begin{equation}
_{\mu ,\nu }^{\gamma ,\sigma }<\varsigma |\overset{\wedge }{n}|\varsigma
>_{\mu ,\nu }^{\gamma ,\sigma }=_{\mu ,\nu }^{\gamma ,\sigma }<\varsigma
|a^{+}a|\varsigma >_{\mu ,\nu }^{\gamma ,\sigma }=\sum\limits_{n=0}^{\infty
}nP_{\mu ,\nu }^{\gamma ,\sigma }(n,\varsigma )=\gamma \Gamma (\nu )\frac{%
|\varsigma |^{2\sigma }}{\Gamma (\mu +\nu )}.  \label{eq7.7b}
\end{equation}

The second order moment of the number of photons in the quantum state $%
|\varsigma >_{\mu ,\nu }^{\gamma ,\sigma }$ is

\begin{equation}
_{\mu ,\nu }^{\gamma ,\sigma }<\varsigma |(\overset{\wedge }{n}%
)^{2}|\varsigma >_{\mu ,\nu }^{\gamma ,\sigma }=_{\mu ,\nu }^{\gamma ,\sigma
}<\varsigma |(a^{+}a)^{2}|\varsigma >_{\mu ,\nu }^{\gamma ,\sigma }=
\label{eq7.7c}
\end{equation}%
\begin{equation*}
\sum\limits_{n=0}^{\infty }n^{2}P_{\mu ,\nu }^{\gamma ,\sigma }(n,\varsigma
)=\frac{(\gamma )_{2}\Gamma (\nu )|\varsigma |^{4\sigma }}{\Gamma (2\mu +\nu
)}+\frac{\gamma \Gamma (\nu )|\varsigma |^{2\sigma }}{\Gamma (\mu +\nu )}.
\end{equation*}

For one-mode quantum field these equations allow us to introduce the
Mandel's parameter \cite{Mandel} $Q_{\mu ,\nu }^{\gamma ,\sigma }$,

\begin{equation}
Q_{\mu ,\nu }^{\gamma ,\sigma }=\frac{_{\mu ,\nu }^{\gamma ,\sigma
}<\varsigma |(\overset{\wedge }{n})^{2}|\varsigma >_{\mu ,\nu }^{\gamma
,\sigma }-(_{\mu ,\nu }^{\gamma ,\sigma }<\varsigma |\overset{\wedge }{n}%
|\varsigma >_{\mu ,\nu }^{\gamma ,\sigma })^{2}}{_{\mu ,\nu }^{\gamma
,\sigma }<\varsigma |\overset{\wedge }{n}|\varsigma >_{\mu ,\nu }^{\gamma
,\sigma }}-1.  \label{eq.7.8}
\end{equation}

Taking into account Eqs.(\ref{eq7.7b}) and (\ref{eq7.7c}) we obtain

\begin{equation}
Q_{\mu ,\nu }^{\gamma ,\sigma }=\left\{ \frac{(\gamma )_{2}\Gamma (\mu +\nu )%
}{\gamma \Gamma (2\mu +\nu )}-\frac{\gamma \Gamma (\nu )}{\Gamma (\mu +\nu )}%
\right\} |\varsigma |^{2\sigma }.  \label{eq7.8a}
\end{equation}

It is easy to see that for $\mu =\nu =\gamma =1$ the Mandel's parameter is

\begin{equation}
Q_{1,1}^{1,\sigma }=0.  \label{eq7.8b}
\end{equation}

The fractional generalized Bell polynomials come into the play in quantum
mechanics when we calculate the diagonal matrix element of the $n^{\mathrm{th%
}}$ power of the number operator $(a^{+}a)^{n}$. Indeed, we have

\begin{equation}
_{\mu ,\nu }^{\gamma ,\sigma }<\varsigma |(a^{+}a)^{n}|\varsigma >_{\mu ,\nu
}^{\gamma ,\sigma }=B_{\mu ,\nu }^{\gamma }(|\varsigma |^{2\sigma },n),
\label{eq7.9}
\end{equation}

where $B_{\mu ,\nu }^{\gamma }(|\varsigma |^{2\sigma },n)$ is defined by Eq.(%
\ref{eq5.1a}) or by Eq.(\ref{eq5.7}).

The Hamiltonian operator $\overset{\wedge }{H}$ of a single-mode quantized
field has the form

\begin{equation*}
\overset{\wedge }{H}=\hbar \omega \overset{\wedge }{n}=\hbar \omega a^{+}a,
\end{equation*}

where $\hbar $ is Planck's constant. Then considering the quantum mechanical
evolution operator $\exp (-i\overset{\wedge }{H}t/\hbar )$

\begin{equation*}
\exp (-i\overset{\wedge }{H}t/\hbar )=\exp (-i\omega ta^{+}a),
\end{equation*}

we see that the diagonal matrix element $_{\mu ,\nu }^{\gamma ,\sigma
}<\varsigma |\exp (-i\overset{\wedge }{H}t/\hbar )|\varsigma >_{\mu ,\nu
}^{\gamma ,\sigma }$ can be written as

\begin{equation}
_{\mu ,\nu }^{\gamma ,\sigma }<\varsigma |\exp (-i\omega ta^{+}a)|\varsigma
>_{\mu ,\nu }^{\gamma ,\sigma }=\sum\limits_{n=0}^{\infty }\frac{1}{n!}%
(-i\omega t)^{n}\ _{\mu ,\nu }^{\gamma ,\sigma }<\varsigma
|(a^{+}a)^{n}|\varsigma >_{\mu ,\nu }^{\gamma ,\sigma }=  \label{eq7.11}
\end{equation}

\begin{equation*}
\sum\limits_{n=0}^{\infty }\frac{1}{n!}(-i\omega t)^{n}B_{\mu ,\nu }^{\gamma
}(|\varsigma |^{2\sigma },n).
\end{equation*}

By comparing this result with Eqs.(\ref{eq6.1}) and (\ref{eq6.3}) we
conclude that

\begin{equation}
_{\mu ,\nu }^{\gamma ,\sigma }<\varsigma |\exp (-i\omega ta^{+}a)|\varsigma
>_{\mu ,\nu }^{\gamma ,\sigma }=E_{\mu ,\nu }^{\gamma }(|\varsigma
|^{2\sigma }(\exp (-i\omega t)-1)).  \label{eq7.12}
\end{equation}

In other words, the diagonal matrix element of the evolution operator $\exp
(-i\omega ta^{+}a)$ in the basis of the coherent states $|\varsigma >_{\mu
,\nu }^{\gamma ,\sigma }$ is expressed as the generating function of the
fractional generalized Bell polynomials. It follows from Eq.(\ref{eq7.12})
that for the special case when \TEXTsymbol{\vert}$\varsigma $\TEXTsymbol{%
\vert}$=1$ the diagonal matrix element of the $n^{\text{th}}$ power of the
number operator $(a^{+}a)^{n}$ yields the fractional generalized Bell number 
$B_{\mu ,\nu }^{\gamma }(n)$,

\begin{equation}
_{\mu ,\nu }^{\gamma ,\sigma }<\varsigma |(a^{+}a)^{n}|\varsigma >_{\mu ,\nu
}^{\gamma ,\sigma }\left\vert _{|\varsigma |=1}\right. =B_{\mu ,\nu
}^{\gamma }(|\varsigma |^{2\sigma },n)\left\vert _{|\varsigma |=1}\right.
=B_{\mu ,\nu }^{\gamma }(n).  \label{eq7.13}
\end{equation}

When $0<\mu \leq 1$, $\nu =\gamma =1$ and $\sigma =\mu $ the Eq.(\ref{eq7.12}%
) gives us

\begin{equation}
_{\mu ,1}^{1,\mu }<\varsigma |\exp ((-i\omega t)a^{+}a)|\varsigma >_{\mu
,1}^{1,\mu }=E_{\mu }(|\varsigma |^{2\mu }(\exp (-i\omega t)-1)),
\label{eq7.14}
\end{equation}

which shows that the diagonal coherent state $|\varsigma >_{\mu ,1}^{1,\mu }$
matrix element of the operator $\exp \left\{ (-i\omega t/\hbar
)a^{+}a\right\} $ is the generating function of the fractional Bell
polynomials \cite{Laskin2}.

When $\mu =\nu =\gamma =1$ the Eq.(\ref{eq7.12}) gives us

\begin{equation}
_{\sigma }<\varsigma |\exp ((-i\omega t)a^{+}a))|\varsigma >_{\sigma }=\exp
(|\varsigma |^{\sigma }(\exp (-i\omega t)-1)),  \label{eq7.15}
\end{equation}

where we introduced the following notation \TEXTsymbol{\vert}$\varsigma
>_{\sigma }=$\TEXTsymbol{\vert}$\varsigma >_{1,1}^{1,\sigma }$.

When $\sigma =1$ the Eq.(\ref{eq7.15}) gives us

\begin{equation}
<\varsigma |\exp ((-i\omega t)a^{+}a)|\varsigma >=\exp (|\varsigma
|^{2}(\exp (-i\omega t)-1)),  \label{eq7.15a}
\end{equation}

where we introduced the following notation \TEXTsymbol{\vert}$\varsigma >=$%
\TEXTsymbol{\vert}$\varsigma >_{1,1}^{1,\sigma }|_{\sigma =1}$.

When \TEXTsymbol{\vert}$\varsigma |=1$ Eqs.(\ref{eq7.15}) and (\ref{eq7.15a}%
) become a relation between the diagonal matrix element of the $n^{\text{th}}
$ power of the number operator $(a^{+}a)^{n}$ in the basis of standard
coherent states $|\varsigma >$\ and Bell numbers $B(n)$,

\begin{equation}
<\varsigma |(a^{+}a)^{n}|\varsigma >|_{|\varsigma |=1}=B(|\varsigma
|^{2},n)|_{|\varsigma |=1}=B(n).  \label{eq7.16}
\end{equation}

The equation (\ref{eq7.16}) was originally obtained in \cite{Katriel}.

\subsubsection{Stretched quantum coherent states}

Let us introduce \textit{stretched quantum coherent states} $|\varsigma
>_{\sigma }$, which are a special case of the states $|\varsigma >_{\mu ,\nu
}^{\gamma ,\sigma }$ when $\mu =\nu =\gamma =1$,

\begin{equation}
|\varsigma >_{\sigma }=|\varsigma >_{1,1}^{1,\sigma }=\exp (-\frac{%
|\varsigma |^{2\sigma }}{2})\sum\limits_{n=0}^{\infty }\frac{\varsigma
^{\sigma n}}{\sqrt{n!}}|n>,\qquad 0<\sigma \leq 1,  \label{eq7.17}
\end{equation}

and the adjoint states $_{\sigma }<\varsigma |$%
\begin{equation}
_{\sigma }<\varsigma |=_{1,1}^{1,\sigma }<\varsigma |=\exp (-\frac{%
|\varsigma |^{2\sigma }}{2})\sum\limits_{n=0}^{\infty }\frac{(\varsigma
^{\sigma \ast })^{n}}{\sqrt{n!}}<n|,\qquad 0<\sigma \leq 1.  \label{eq7.18}
\end{equation}

It is easy to see that stretched quantum coherent state $|\varsigma
>_{\sigma }$ is eigenstate of photon field annihilation operator $a$. \
Indeed, we have

\begin{equation}
a|\varsigma >_{\sigma }=\exp (-\frac{|\varsigma |^{2\sigma }}{2}%
)\sum\limits_{n=0}^{\infty }\frac{\varsigma ^{\sigma n}}{\sqrt{n!}}\sqrt{n}%
|n-1>=  \label{eq7.20}
\end{equation}

\begin{equation*}
\exp (-\frac{|\varsigma |^{2\sigma }}{2})\sum\limits_{m=0}^{\infty }\frac{%
\varsigma ^{\sigma (n+1)}}{\sqrt{n!}}|n>=\varsigma ^{\sigma }|\varsigma
>_{\sigma },
\end{equation*}

that is

\begin{equation}
a|\varsigma >_{\sigma }=\varsigma ^{\sigma }|\varsigma >_{\sigma },
\label{eq.7.20a}
\end{equation}

which shows that the eigenvalue of the photon field annihilation operator $a$
acting on the stretched quantum coherent state $|\varsigma >_{\sigma }$ is $%
\varsigma ^{\sigma }$.

Using Eq.(\ref{eq7.1}) and the following representation for the eigenvector $%
|n>$ of the photon number operator

\begin{equation}
|n>=\frac{(a^{+})^{n}}{\sqrt{n!}}|0>,  \label{eq7.21}
\end{equation}

where $|0>$ is the vacuum state with $\varsigma =0$, we obtain the
alternative expression for the stretched quantum coherent states

\begin{equation}
|\varsigma >_{\sigma }=\exp (-\frac{|\varsigma |^{2\sigma }}{2}%
)\sum\limits_{n=0}^{\infty }\frac{(\varsigma ^{\sigma }a^{+})^{n}}{\sqrt{n!}}%
|0>,\qquad 0<\sigma \leq 1.  \label{eq7.22}
\end{equation}

The \textit{overcompleteness} relation has the form

\begin{equation}
_{\sigma }<\eta |\varsigma >_{\sigma }=\exp \{-\frac{|\eta |^{2\sigma }}{2}-%
\frac{|\varsigma |^{2\sigma }}{2}+\eta ^{\sigma \ast }\varsigma ^{\sigma }\}.
\label{eq7.23}
\end{equation}

The squared modulus $|<n|\varsigma >_{\sigma }|^{2}$ of projection of the
stretched coherent state $|\varsigma >_{\sigma }$ onto the number state $|n>$
gives us the probability $P_{\sigma }(n,\varsigma )$ that $n$ photons will
be found in the coherent state $|\varsigma >_{\sigma },$

\begin{equation}
P_{\sigma }(n,\varsigma )=|<n|\varsigma >_{\sigma }|^{2}=\frac{|\varsigma
|^{2\sigma n}}{n!}\exp (-|\varsigma |^{2\sigma }).  \label{eq7.19}
\end{equation}

Therefore, the stretched quantum coherent states $|\varsigma >_{\sigma }$
are normalized due to the normalization condition for $P_{\sigma
}(n,\varsigma )$ 

\begin{equation*}
\sum\limits_{n=0}^{\infty }|<n|\varsigma >_{\sigma
}|^{2}=\sum\limits_{n=0}^{\infty }P_{\sigma }(n,\varsigma )=1.
\end{equation*}

The mean number of photons in the quantum state $|\varsigma >_{\sigma }$ is

\begin{equation*}
_{\sigma }<\varsigma |\overset{\wedge }{n}|\varsigma >_{\sigma =\sigma
}=<\varsigma |a^{+}a|\varsigma >_{\sigma }=\sum\limits_{n=0}^{\infty
}nP_{\sigma }(n,\varsigma )=|\varsigma |^{2\sigma }.
\end{equation*}

Taking into account Eq.(\ref{eq7.8b}) we conclude that stretched quantum
coherent states obey the Poisson statistics.

\subsubsection{Stretched displacement operator}

We introduce the \textit{stretched displacement operator} $D_{\sigma
}(\varsigma )$ as follows

\begin{equation}
D_{\sigma }(\varsigma )=\exp \{\varsigma ^{\sigma }a^{+}-\varsigma ^{\sigma
\ast }a\},  \label{eq7.24}
\end{equation}

where $\varsigma $ is a complex number, $a^{+}$ and $a$ are photon field
creation and annihilation operators.

When $\sigma =1$ the operator $D_{\sigma }(\varsigma )|_{\sigma =1}$ becomes
the well-known displacement operator $D(\varsigma )$ for the Schr\"{o}%
dinger-Glauber coherent states \cite{Schrodinger} - \cite{Klauder},

\begin{equation*}
D(\varsigma )=D_{1}(\varsigma )=\exp \{\varsigma a^{+}-\varsigma ^{\ast }a\}.
\end{equation*}

The stretched displacement operator $D_{\sigma }(\varsigma )$ is an unitary
operator, i.e.

\begin{equation}
D_{\sigma }^{+}(\varsigma )D_{\sigma }(\varsigma )=1,  \label{eq7.24a}
\end{equation}

where the sign "$+"$ stands for Hermitian conjugation of an operator. Using
the Baker-Campbell-Hausdorff formulae for operators $\widehat{A}$ and $%
\widehat{B}$

\begin{equation*}
e^{\widehat{A}+\widehat{B}}=e^{\widehat{A}}e^{\widehat{B}}e^{-\frac{1}{2}[%
\widehat{A},\widehat{B}]},
\end{equation*}

such that%
\begin{equation*}
\lbrack \widehat{A},[\widehat{A},\widehat{B}]]=[\widehat{B},[\widehat{A},%
\widehat{B}]]=0,
\end{equation*}

we come to alternative representations for the stretched displacement
operator $D_{\sigma }(\varsigma )$

\begin{equation}
D_{\sigma }(\varsigma )=\exp (-\frac{|\varsigma |^{2\sigma }}{2})\exp
(\varsigma ^{\sigma }a^{+})\exp (-\varsigma ^{\sigma \ast }a).
\label{eq7.25a}
\end{equation}

or

\begin{equation}
D_{\sigma }(\varsigma )=\exp (\frac{|\varsigma |^{2\sigma }}{2})\exp
(-\varsigma ^{\sigma \ast }a)\exp (\varsigma ^{\sigma }a^{+}).
\label{eq7.25b}
\end{equation}

To show that the action of the stretched displacement operator on vacuum
state gives the stretched coherent state $|\varsigma >$

\begin{equation}
|\varsigma >_{\sigma }=D_{\sigma }(\varsigma )|0>,  \label{eq7.26}
\end{equation}

we perform the following chain of transformations

\begin{equation*}
D_{\sigma }(\varsigma )|0>=\exp (-\frac{|\varsigma |^{2\sigma }}{2})\exp
(\varsigma ^{\sigma }a^{+})\exp (-\varsigma ^{\sigma \ast }a)|0>=
\end{equation*}

\begin{equation*}
\exp (-\frac{|\varsigma |^{2\sigma }}{2})\exp (\varsigma ^{\sigma
}a^{+})|0>=\exp (-\frac{|\varsigma |^{2\sigma }}{2})\sum\limits_{n=0}^{%
\infty }\frac{\varsigma ^{\sigma n}(a^{+})^{n}}{n!}|0>=
\end{equation*}

\begin{equation*}
\exp (-\frac{|\varsigma |^{2\sigma }}{2})\sum\limits_{n=0}^{\infty }\frac{%
\varsigma ^{\sigma n}}{\sqrt{n!}}|n>=|\varsigma >_{\sigma },
\end{equation*}

where Eq.(\ref{eq7.21}) was used.

It is easy to see that the following equations hold

\begin{equation}
\lbrack a,D_{\sigma }(\varsigma )]=\varsigma ^{\sigma }D_{\sigma }(\varsigma
),  \label{eq7.27}
\end{equation}

and

\begin{equation}
D_{\sigma }^{+}(\varsigma )aD_{\sigma }(\varsigma )=a+\varsigma ^{\sigma
},\qquad \qquad D_{\sigma }^{+}(\varsigma )a^{+}D_{\sigma }(\varsigma
)=a^{+}+\varsigma ^{\sigma \ast }.  \label{eq7.28}
\end{equation}

For $\sigma =1$, all the above equations turn into the well-known equations
for Schr\"{o}dinger-Glauber coherent states of the standard quantum
mechanics \cite{Schrodinger} - \cite{Klauder}.

\section{Conclusion}

We introduce and develop a new family of fractional discrete stochastic
processes and probability distribution functions based on a three-parameter
generalized Mittag-Leffler function known as the Prabhakar function. The new
family includes the famous Poisson probability distribution and the
fractional Poisson probability distribution. Fundamentals and applications
of the new fractional generalized counting probability distribution were
developed. As applications we built a fractional generalized compound
process, fractional generalized Bell polynomials and fractional generalized
Bell numbers, fractional generalized Stirling numbers, and obtained new
equations for probability distribution of waiting time.

We introduce stretched quantum coherent states and show that they support
fundamental properties of the celebrated Schr\"{o}dinger-Glauber coherent
states.

Under certain sets of fractality parameters, all our new results are
transformed into the well-known equations for the Poisson process and
probability distribution, the fractional Poisson process and probability
distribution, Bell polynomials and numbers, fractional Bell polynomials and
numbers, Stirling and fractional Stirling numbers of the second kind, and
well-known quantum coherent states.

\end{document}